\documentclass[12pt]{article}
\usepackage{amsmath, amssymb,amsfonts}
\usepackage{graphicx}
\usepackage{amsmath}
\usepackage{amsfonts}
\usepackage{amssymb}

\newtheorem{teo}{Theorem}

\newtheorem{lemma}{Lemma}

\newtheorem{definition}{Definition}
\newtheorem{remark}{Remark}

\begin{document}

\title{Flux of superconducting vortices through a domain}
\author{S.N. Antontsev and N.V. Chemetov \\
{\small CMAF/Universidade de Lisboa,  Portugal and }\\
{\small University of Sao Paulo, 14040-901 Ribeirao Preto - SP, Brazil } }

\date{}
\maketitle
\tableofcontents

\textbf{Abstract.} {\small We study a mean-field model of superconducting
vortices in a II-type superconductor. We prove the global solvability of the
problem about a flux of superconducting vortices through a given domain and
also analyse the regularity of weak solutions. }



\section{The statement of the problem}

\label{sec1}

Certain materials, when cooled through a critical temperature, exhibit a
superconducting state in which they have the ability to conduct electric
currents without resistance. Many models have been proposed to describe the
behaviour of superconductors: the microscopic theory of Bardeen and
Cooper-Schreiffer, the mesoscopic theories of London and Ginzburg-Landau,
the macroscopic critical state theories, like the Bean model. For the
description and relationship of these models we refer to the excellent works
of S.J. Chapman \cite{Chap00-1}, \cite{Chap00-2}. In this paper we
investigate a mean field model of superconducting vortices in a II-type
superconductor, formulated in \cite{Chap95}. In 1-D and 2-D cases the
following system describes the evolution of the vorticity $\omega$ and the
average magnetic field $h$ of the superconducting sample:
\begin{eqnarray}
\omega_{t}+\mbox{div}\,(\omega\vec{v})=0, & (\vec{x},t)\in\Omega_{T}:=\Omega%
\times(0,T),  \label{eq1} \\
-\Delta h+h=\omega, & (\vec{x},t)\in\Omega_{T},  \label{eq2} \\
\vec{v} = -\nabla h, & (\vec{x},t)\in\Omega_{T},  \label{eq3}
\end{eqnarray}
with the following boundary conditions
\begin{eqnarray}
\vec{v}\,\vec{n} = a(\vec{x},t), & (\vec{x},t)\in\Gamma_{T}:=\Gamma%
\times(0,T),  \label{eq4} \\
\omega= b(\vec{x},t), & (\vec{x},t)\in\Gamma^{-}_{T}:=\Gamma^{-}\times(0,T)
\label{eq5}
\end{eqnarray}
and the initial conditions
\begin{eqnarray}
\omega(\vec{x},0)=\omega_{0}(\vec{x}),\quad\vec{x}\in\Omega.  \label{eq6}
\end{eqnarray}
Here $\vec n$ denotes the outward normal to the $C^{2+\gamma }$-smooth
boundary $\Gamma=\Gamma^{+}\cup \Gamma^{0}\cup \Gamma^{-}$ of the bounded domain $%
\Omega\subset\mathbb{R}^{n},$ with some $\gamma>0$; $\Gamma^{+}$ is the part
of the boundary $\Gamma$, where $\vec{v}\,\vec{n}=a >0$; $\Gamma^{0}$ is the
part of $\Gamma$, where $\vec{v}\,\vec{n}=a= 0$; $\Gamma^{-}$ is the
part of $\Gamma$, where $\vec{v}\,\vec{n}=a < 0$ and $meas(\Gamma^{+})\ne0
,\,\,\, meas(\Gamma^{-}) \ne0 $. Let us note that the physical statement of
the system is valid just for $n=1$ and $n=2$, but in this article we shall
consider the common case, when $n\geqslant 1$. Discussions of different
types of boundary and initial conditions can be found in the above mentioned
papers \cite{Chap00-1}, \cite{Chap00-2} and \cite{Chap95}. Note that the
system (\ref{eq1})-(\ref{eq3}) is hyperbolic-elliptic type for unknown
functions $\{\omega,h\}.$

Many articles on (\ref{eq1})-(\ref{eq3}) have been written in the last ten
years. The problem has been studied through various approaches ranging from
asymptotic analysis, numerical simulations and rigorous mathematical
analysis. As we know the first positive existence result for the system (\ref%
{eq1} )-(\ref{eq3}), (\ref{eq6}) was obtained by Huang and Svobodny \cite%
{HuSvo} in the case of the Cauchy problem for the domain $\Omega=\mathbb{R}%
^{2}$, using the method of characteristics and potentials. In the articles
\cite{SchatSty99}-\cite{Bri}, the system (\ref{eq1})-(\ref{eq3}), (\ref{eq6}%
) was considered in the case of a bounded domain $\Omega$ with boundary
conditions
\begin{eqnarray}
h(\vec{x},t) = a, & (\vec{x},t)\in\Gamma_{T},  \notag \\
\omega(\vec{x},t)=b, & (\vec{x},t)\in\Gamma_{T}^-,  \label{eq4555}
\end{eqnarray}
where $a, \, b$ are given \textit{constants}. In \cite{SchatSty99} the
existence of a weak solution of (\ref{eq1})-(\ref{eq3}), (\ref{eq6})-(\ref%
{eq4555}) for $b=0$ was deduced by an approach based on a parabolic-elliptic
approximation. In 1-D case the authors have shown the uniqueness of the
solution. Also in this work and in \cite{Cha-R-S} the existence of a steady
state solution $\{\omega, h\}$ of (\ref{eq1})-(\ref{eq3}), (\ref{eq6})-(\ref%
{eq4555}) was demonstrated. In \cite{ElSchatSto00} the solvability of (\ref%
{eq1})-(\ref{eq3}), (\ref{eq6})-(\ref{eq4555}) with an additional curvature
term in (\ref{eq1}) was shown. In \cite{ElSty00}, the case of flux pinning
and boundary nucleation of vorticity (i.e. the \textit{constants} $a, b>0$)
was considered, where existence and uniqueness of a solution in 1-D case
were proved. The discretization, the convergence of discretizated solutions
and numerical simulations of (\ref{eq1})-(\ref{eq3}) were extensively
studied in articles \cite{X.Chen}, \cite{C.M.Elliott-V.M.Styles2}, \cite%
{Styles}, \cite{Bri}. In all mentioned works the boundary conditions (\ref%
{eq4555}) were taken as \textit{zero or constant } valued. But the evolution
of vortices arises under penetration of the magnetic field into
superconducting bodies, i.e. the case of non zero or non constant boundary
conditions is really physical.

\textbf{Remark}: This is a version of the article entitled ” Flux of superconducting vortices through a domain” published in SIAM Journal on Mathematical Analysis, 39 (2007), pp. 263-280, we refer to \cite{AC1}.  Since the publication of this paper, other articles on this subject have appeared, such as Antontsev, Chemetov \cite{AC2}. Moreover in the articles  \cite{chem1, chem2}  a significant progress in the study of the superconductive model has been done, where the Kruzkov approach  has been generalized  using the kinetic method. We also mention the articles \cite{chem3}-\cite{chem5}, where the last kinetic model has been applied to porous media models. 

Since the question of the trace value of integrable solution is important in the study of the problem we also refer to  \cite{C}-\cite{CC6}. The stochastic case for  system similar to the one treated here \cite{ACC21}, \cite{CC0}.

In the present article we consider the problem (\ref{eq1})-(\ref{eq6}) for
the general case, when $a(\vec x, t), \, b(\vec x, t) $ are \textit{arbitrary%
} given functions of the conditions (\ref{eq4})-(\ref{eq5}). We prove the
existence of solution $\{\omega,h\}$ with a natural smooth restriction on
the data of the problem.

\textbf{Terminology and notations.} In accordance with the notations,
introduced in the books \cite{LadySolonUral68}, \cite{LadyUral68}, we shall
use the Sobolev's spaces $L_{q} (\Omega),$ $W^{l}_{q} (\Omega)$, $%
W^{l,\,m}_{q} (\Omega_{T} )$, $W^{l}_{q}(\Gamma )$, $\,\,\,W^{l,\,m}_{q}(%
\Gamma_{T} )\,\,$ for $\,\,l, q \geqslant1$, $m\geqslant 0$ and the H\"{o}%
lder spaces $C^l (\Omega)$, $C^{l,\,m} (\Omega_{T} )$, $C^{l,\,m}
(\Gamma_{T} )\,\,$ for $\,\,l, m \geqslant 0$, where the index $l$
corresponds to the variable $\vec x$ and $m$ to the variable $t$ ($\,\,l,q,m$
are integer or non-integer).

Let $B$ be a given Banach space. We denote by $C(0,T; B)$ the space of
continuous functions from $[0,T]$ into the Banach space $B$ with the norm
\begin{eqnarray*}
||u||_{C(0,T; \,\,B)} =\max_{t \in [0,T]} ||u(t)||_B
\end{eqnarray*}
and by $L_p (0,T; B)$ for any $1\leqslant p \leqslant \infty$, the space of
measurable functions from $[0,T]$ into the Banach space $B$ with $p$-th
power summable on $[0,T]$, with the norm
\begin{eqnarray*}
||u||_{L_p (0,T;\,\, B)} =(\int _0^T ||u(t)||_B ^{p} d t )^{\frac{1}{p}%
}\;\;\;\;\; ( =ess \sup_{ 0< t < T } ||u(t)||_B \mbox{ if } p=\infty ).
\end{eqnarray*}

\textbf{Regularity of data.} \textit{We assume that the datum $\,\,\,a\,\,\,$
satisfies the condition
\begin{equation}
a\in C^{\eta ,\, \theta } (\Gamma_{T}) \,\, \mbox{ for some } \eta , \,
\theta \in (0,1)  \label{eq00sec1}
\end{equation}
and the data $b ,\,\, \omega_{0}$ admit an extension $\breve\omega$, defined
on all domain $\overline\Omega_{T} $, such that
\begin{eqnarray}
\breve\omega(\vec x,t )=b(\vec x ,t ), \,\,\,\, (\vec x, t
)\in\Gamma_{T}^{-} ;\;\;\;\;\,\,\, \breve\omega(\vec x , 0 ) = \omega_{0}
(\vec x ) ,\,\,\,\, \vec x \in\Omega;  \label{eq01sec1}
\end{eqnarray}
and $\breve\omega$ satisfies
\begin{eqnarray}
0 \leqslant\breve\omega(\vec x, t) \leqslant\aleph\,\,\,\,\, &
\mbox{ for
a.e. }& (\vec x, t) \in\overline\Omega_{T} ,\;\;\;  \label{eq02sec1} \\
( || \partial_{t} \breve\omega(\cdot, t) ||_{L_{1} (\Omega)}+
||\bigtriangledown\breve\omega(\cdot, t) ||_{L_{\infty}(\Omega)}) &\in L_{1}
(0,T) .& \,\,\,  \label{eq03sec1}
\end{eqnarray}
where $\aleph$ is a constant. }

\begin{remark}
The assumption (\ref{eq02sec1}) implies that the functions $\omega_{0}$ and $%
b$ are positive and bounded above by the constant $\aleph$ on $\overline
\Omega$ and $\Gamma^{-}_{T}$, respectively. In section \ref{sec6} we give
sufficient conditions on the data $b,\,\, \omega_{0}$, that permit the
existence of an extension $\breve \omega$ on all domain $\overline\Omega_{T}$%
, satisfying the conditions (\ref{eq01sec1})-(\ref{eq03sec1}).
\end{remark}

Let us give the definition of the weak solution of our problem.

\begin{definition}
\label{def1sec1} A pair of functions $\{\omega,h\}$ is said to be a weak
solution of the problem (\ref{eq1}) -- (\ref{eq6}), if $\omega\in L_{\infty
}(\Omega_{T}), \; $ $h\in L_{\infty}(0,T;\; C^{1+\alpha} (\overline\Omega))
\, \cap\, L_{2} (0,T;\; W^{2}_{2}(\Omega^{\prime})) $ for some $\alpha \in
(0,1)$ and for any sub domain $\Omega^{\prime}$ of $\Omega $ such that $%
\overline{ \Omega^{\prime}} \subset\Omega$ and the following equalities
\begin{eqnarray}
\!\!\!\int_{\Omega_{T}}\omega(\psi_{t} + \vec{v}\,\nabla\psi)\; d\vec{x} dt
& + \int_{\Omega}\omega_{0}\; \psi(\vec x,0) \; d\vec{x} = -\int
_{\Gamma^{-}_{T}} a\;b\;\psi\;d\vec{x} dt,  \label{eq7} \\
-\Delta h+h = \omega & \mbox{a.e. in } \Omega_{T},  \label{eq8} \\
\vec{v} =-\nabla h & \mbox{a.e. in } \Omega_{T},  \label{eq9} \\
-\nabla h\,\vec{n}=a & \mbox{a.e. on } \Gamma_{T},  \label{eq10}
\end{eqnarray}
hold for arbitrary function $\psi\in H^{1}(\Omega_{T}),$ such that
\begin{eqnarray}
\psi(\vec{x},T)=0 \mbox{ for } \vec{x}\in\Omega \;\;\;\;\mbox{ and }\;\;\;\;
\psi(\vec{x},t)=0 \mbox{ for } (\vec {x},t)\in\Gamma^{+}_{T}\cup \Gamma^{0}_{T}.  \label{psi}
\end{eqnarray}
\end{definition}

Our main result in this work is the following theorem.

\begin{teo}
\label{teo2sec1} If the data $a,\,\, b ,\,\, \omega_{0}$ satisfy (\ref%
{eq00sec1})-(\ref{eq03sec1}), then there exists at least one weak solution $%
\{\omega,h\}$ of the problem (\ref{eq1}) -- (\ref{eq6}). Moreover, we have
\begin{eqnarray}
\omega\in L_\infty (\Omega_T),\; \partial_t\omega\in L_\infty
(0,T;H^{-1}(\Omega))  \label{A2}
\end{eqnarray}
and
\begin{eqnarray}
h&\in& L_\infty (0,T;\, C^{1+\eta} (\overline{\Omega} )) \cap L_\infty
(0,T;\,W^2_q(\overline{\Omega^{\prime }} )),\,  \notag \\
h, \,\, \bigtriangledown h &\in& C^{0 , \, \theta } (\overline{%
\Omega^{\prime }}\times [0,T] )  \label{A22}
\end{eqnarray}
for any $q\in(1,\infty)$ and for any $%
\Omega^{\prime },$ such that $\overline{\Omega^{\prime }}\subset \Omega$.
\end{teo}

\begin{remark}
The similar result is valid for both Dirichlet  and Robin boundary conditions 
(in accordance with "www.wikipedia.org"
 terminology), instead of  (\ref{eq10}), on the function $h$. For simplicity, we prove
this theorem for the dimension $n\geqslant 2$.
\end{remark}

We divide the proof of this theorem into 4 steps:

$1^{st}$ step: in section \ref{sec2} we remember well-known results from the
theory of elliptic equations, that will be very useful in the sequel.

$2^{nd}$ step: in section \ref{sec3}, we show that an approximate problem to
(\ref{eq1}) - (\ref{eq6}), which is of parabolic-elliptic type, is solvable,
if the initial and boundary conditions are smooth; the reasoning is based on
the Schauder fixed-point theorem.

$3^{d}$ step: in section \ref{sec4}, by using the principle of
maximum-minimum and the results of potential theory for elliptic equations,
mentioned in the section \ref{sec2}, we deduce $L_{\infty}$ - estimates and
some additional a priori estimates for the solutions of the approximate
problem of (\ref{eq1})-(\ref{eq6}).

$4^{th}$ step: in section \ref{sec5} we use a compactness argument to
establish the existence of a solution for the original problem.

\section{Auxiliary assertions. The elliptic boundary value problem}

\label{sec2} \setcounter{equation}{0}

In this section for the convenience of the reader we give well-known results
on elliptic equations.

According to the potential theory for elliptic equations \cite{Vlad00}, the
solution of the problem
\begin{eqnarray}
\left\{
\begin{array}{lll}
-\Delta h_{1}+h_{1}=F, & \vec x\in\Omega; &  \\
&  &  \\
-\nabla h_{1}\,\vec n=0, & \vec x\in\Gamma &
\end{array}
\right.  \label{eq13''sec3}
\end{eqnarray}
can be written in the form $h_{1}(\vec{x})=(K_{1} \ast F) (\vec{x}%
):=\int_{\Omega}K_{1}(\vec{x},\vec{y}) \,\,F (\vec y)\; d\vec y $ and the
solution of the problem
\begin{eqnarray}
\left\{
\begin{array}{lll}
-\Delta h_{2}+h_{2}=0, & \vec x\in\Omega; &  \\
&  &  \\
-\nabla h_{2}\,\vec n=g, & \vec x\in\Gamma &
\end{array}
\right.  \label{eq13'sec3}
\end{eqnarray}
can be written in the form $h_{2}=(K_{2} \ast g )(\vec{x}):=\int_{%
\Gamma}K_{2}(\vec{x},\vec{y})g(\vec y)\; d\vec y. $ Since $\Gamma$ is $%
C^{2+\gamma}$-smooth, the kernels $K_{1},\,\, K_{2}$ satisfy the
inequalities
\begin{eqnarray*}
\left| K_{i}(\vec{x},\vec{y})\right| &\leqslant& C|x-y|^{2-n}, \\
\left| \nabla_{\vec x}K_{i}(\vec{x},\vec{y})\right| &\leqslant&
C|x-y|^{1-n},\;\;i=1,2,\quad \mbox{ for any } \vec x, \vec y \in\Omega
\end{eqnarray*}
in $n$-dimensional case with $n>2$ and
\begin{eqnarray*}
\left| K_{i}(\vec{x},\vec{y})\right| &\leqslant& C|ln|x-y||, \\
\left| \nabla_{\vec x }K_{i}(\vec{x},\vec{y})\right| &\leqslant&
C|x-y|^{-1},\;\;i=1,2, \quad \mbox{ for any } \vec x, \vec y \in\Omega
\end{eqnarray*}
in $2$-dimensional case. According to the theory of elliptic equations \cite{LadyUral68}, p. ??;
 to the potential theory
 \cite{Vlad00}, p.??;
 p. 191, Lemma 1.4 of \cite{kato} and embedding
theorems of Sobolev \cite{LadySolonUral68},
the operators $F \to K_{1} \ast F $, $g \to K_{2} \ast g $ possess the
following properties.

\begin{lemma}
\label{teo2sec2} For any $n$-dimensional case with $n\geqslant 2$:

1) The function $h_{1}=K_{1}\ast F$ satisfies the following estimates
\begin{eqnarray}  
h_{1}&\geqslant& 0 \;\text{ in }\; \Omega,
\;\;\;\;\;\;\;\;\;\;\;\;\;\;\;\;\;\;\;\;\;\;\;\;\;\;\;\; \text{if}\;
F\geqslant 0 \text{ a.e. in } \Omega ;  \label{ell-1} \\
||h_{1}||_{C^{1+\alpha }(\overline{\Omega })}&\leqslant&
C||F||_{L_{p}(\Omega )},\;  \alpha = p \in (0,1) \;\; \mbox{ if }
n<p<\infty;  \label{ell-2} \\
||h_{1}||_{C^{1+\alpha }(\overline{\Omega })}&\leqslant&
C||F||_{L_{\infty}(\Omega )},\;\; \forall \; \alpha \in (0,1); \label{ell-21} \\
||h_{1}||_{W_{p}^{2}(\Omega )}&\leqslant& C||F||_{L_{p}(\Omega
)},\;\;\;\;\;\;\;\;\;\;\;\;\;\;\;\;\;\;\;\;\;\;\; \mbox{ if } 1<p<\infty ; \label{ell-3} \\
||h_{1}||_{W^1_{1}(\Omega )}&\leqslant& C||F||_{L_{1}(\Omega )}. \label{ell-4}
\end{eqnarray}

2) The function $h_{2}=K_{2}\ast g$ satisfies the following estimates
\begin{eqnarray}  
||h_{2}||_{C^{1+\alpha }(\overline{\Omega })}&\leqslant&
C||g||_{C^{\alpha }(\Gamma
)},\;\;\;\;\;\;\;\;\;\;\;\;\;\;\;\;\;\;\;\;\;\;\;\;\;\; \mbox{ if }
\;0<\alpha <1; \label{ell-5} \\
||h_{2}||_{C^{l}(\overline{\Omega }^{\prime })}&\leqslant&
C||g||_{L_{1}(\Gamma )}, \label{ell-6}
\end{eqnarray}
for any $\Omega ^{\prime },\;$ such that $\overline{\Omega }^{\prime
}\subset \Omega $ and any $l\geqslant 0$.
\end{lemma}

The proof of this lemma can be found in \cite{LadyUral68} and \cite{Vlad00}.

\section{Construction of approximate solutions}

\label{sec3} \setcounter{equation}{0}

Let $\breve\omega^{{\varepsilon}}$, $a^{\varepsilon}$ for ${\varepsilon}>0$
be smooth approximations of the functions $\breve\omega$, $a$, such that
\begin{eqnarray}
0 \leqslant&\breve\omega^{\varepsilon}(\vec x, t)&\leqslant\aleph\,\,\,\,\, %
\mbox{ for  } \,\,\,\,\, (\vec x, t) \in\overline\Omega_{T};  \notag \\
0 < &a^{\varepsilon}(\vec x, t)& \,\,\,\,\, \,\,\,\,\, \,\,\,\,\,%
\mbox{ for
} \,\,\,\,\, (\vec x, t) \in \Gamma^+_{T} ;  \notag \\
0=&a^{\varepsilon}(\vec x, t)&  \,\,\,\,\, \,\,\,\,\, \,\,\,\,\,%
\mbox{ for
} \,\,\,\,\, (\vec x, t) \in \Gamma^0_{T} ;  \notag \\
&a^{\varepsilon}(\vec x, t)& < 0 \,\,\,\,\, \mbox{ for  } \,\,\,\,\, (\vec
x, t) \in \Gamma^-_{T}.  \label{eq122sec2}
\end{eqnarray}
We require also that
\begin{eqnarray}
\left\{
\begin{array}{ll}
\breve\omega^{\varepsilon} (\vec{x} , t) \mathop{\longrightarrow}\limits_{{%
\varepsilon}\to 0} \breve\omega(\vec{x}, t ) & \mbox{ a.e. in }
\overline\Omega_{T}; \\
&  \\
||a^{\varepsilon}-a\,||_{C^{\eta , \theta } (\Gamma_{T})} %
\mathop{\longrightarrow}\limits_{{\varepsilon}\to 0} 0. &
\end{array}
\right.  \label{ABsec2}
\end{eqnarray}
and
\begin{eqnarray}
\int_{0}^{T} ( || \partial_{t} (\breve\omega^{\varepsilon} - \breve\omega )
(\cdot, t) ||_{L_{1} (\Omega)}+ ||\bigtriangledown
(\breve\omega^{\varepsilon} - \breve\omega )(\cdot, t)
||_{L_{\infty}(\Omega)})\;dt \mathop{\longrightarrow}\limits_{{\varepsilon}%
\to 0} 0.  \label{eq1122sec2}
\end{eqnarray}%
Let us fix a positive number $R$. Next we construct the pair $\{ \omega_{{{{%
\varepsilon} }},R}, h_{{\varepsilon},R}\}$ as a solution of auxiliary
\textbf{Problem} $\mathbf{P_{{\varepsilon},R}}$. For the sake of simplicity,
in this section we suppress the dependence $\omega_{{\varepsilon},R}, h_{{%
\varepsilon},R}$ on ${\varepsilon}, R$ and write $\omega$ and $h$.

Let us consider the problem which is a coupling of the two following systems.

\textbf{Problem} $\mathbf{P_{{\varepsilon},R}}$ \textit{Find $\omega\in
W^{2,1}_{2}(\Omega_{T})$, satisfying the system
\begin{eqnarray}
\left\{
\begin{array}{lll}
\omega_{t} + \mbox{div} (\omega\vec{v})={\varepsilon}\Delta\omega \;\;\;%
\mbox{ and }\;\;\;\; \vec{v}=-\nabla h \;\;\;\;\;\;\;\;\;\;\mbox{ for }%
\;\;\;(\vec{x},t)\in\Omega_{T}; &  &  \\
&  &  \\
\omega(\vec{x},t)=\breve\omega^{\varepsilon}(\vec{x},t), \; (\vec{x}%
,t)\in\Gamma_{T};\;\;\;\;\; \;\;\;\;\;\omega(\vec{x},0)=\breve\omega^{%
\varepsilon}(\vec{x}, 0), \; \vec{x}\in \Omega &  &
\end{array}
\right.  \label{eq3sec2}
\end{eqnarray}
and find $h\in W^{2}_{2} (\Omega) $, satisfying the system
\begin{eqnarray}
\left\{
\begin{array}{ll}
-\Delta h+h=[\omega]_{R}, & (\vec{x},t)\in\Omega_{T}; \\
&  \\
-\nabla h\,\vec{n}= a^{\varepsilon}(\vec{x},t), & (\vec{x},t)\in \Gamma_{T},%
\end{array}
\right.  \label{eq4sec2}
\end{eqnarray}
where $[\cdot]_{R}$ is the cut-off function defined as $[\phi]_{R}:=\max%
\Big\{0,\;\min\{R,\phi\}\Big\}$. }

To prove the solvability of \textbf{Problem} $\mathbf{P_{{\varepsilon},R}}$
we use the Schauder fixed point argument. Let us introduce the class of
functions
\begin{eqnarray}
{\mathcal{M }}=\{\omega(\vec{x},t)\in C (0,T;\;
L_{2}(\Omega)):||\omega||_{C(0,T;\; L_{2} (\Omega))}\leqslant M \},
\label{eq5sec2}
\end{eqnarray}
where an exact value of $M$ will be determined below. In the first we define
the operator $T_{1}$, that transforms a "fixed" vorticity into the
corresponding superconductive field
\begin{eqnarray}
{\mathcal{M}}\ni\widetilde{\omega}\mapsto T_{1}[\widetilde{\omega}]=h
\label{eq6sec2}
\end{eqnarray}
as the solution of (\ref{eq4sec2}), where instead of $\omega$ we put the
chosen $\widetilde{\omega}$. By (\ref{eq13''sec3}), (\ref{eq13'sec3}) the solution $h$ can
be represented in the form
\begin{equation}
h(\vec{x},t) =(K_{1} \ast [\widetilde{\omega}]_{R}) (\vec{x},t) +(K_{2} \ast
a^{\varepsilon}) (\vec{x},t)  \label{eq66sec2}
\end{equation}
for a.e. $(\vec x, t)\in \Omega_T$ and by (\ref{ell-2}), (\ref{ell-5}) of Lemma \ref{teo2sec2} we derive the estimate
\begin{eqnarray}
||h||_{L_{\infty}(0,T; \,\, C^1 (\overline{\Omega } ))}  \leqslant C ( ||\,[\widetilde{\omega}%
]_{R}\,||_{L_{\infty}(\Omega_{T})}+||a^{{\varepsilon}}||_{L_{\infty}(0,T;\;
C^{\eta }(\Gamma))})  \notag \\
\leqslant C ( R +1) .&  \label{eq14sec2}
\end{eqnarray}
The second operator $T_{2}$ describes the evolution of the vorticity
\begin{eqnarray}
\vec{v}=-\nabla h\mapsto T_{2} [h]=\omega,  \label{eq8sec2}
\end{eqnarray}
where $\omega$ is the solution of (\ref{eq3sec2}) for the found
''superconductive field'' $h$ in (\ref{eq6sec2}). Taking into account \ref{eq14sec2} and  the results of \cite%
{LadySolonUral68}, Theorem 4.1, p.153 and Theorem 4.2, p.160, there exists
an unique weak solution $\omega$ of (\ref{eq3sec2}) such that
\begin{eqnarray}
||\omega||_{W^{1, 0}_{2} (\Omega_{T})} &\leqslant& C(R),  \label{eq9sec2} \\
||\omega(\cdot,t_{1})-\omega(\cdot,t_{2})||_{L_{2}(\Omega)} &\leqslant&
\phi_{R} (|t_{1}-t_{2}|),\quad\forall t_{1}, t_{2}\in[0,T],  \label{eq10sec2}
\end{eqnarray}
where the constant $C(R)$ and the function $\phi_{R}=\phi_{R}(t)\in
C([0,T]),\; \phi(0)=0$ depend only on $R$ and the data $\breve\omega ^{{%
\varepsilon} },\,\, a^{\varepsilon}$. Setting $M:=C(R)$ in (\ref{eq5sec2})
and $T=T_2\circ T_1$, by (\ref{eq5sec2}), (\ref{eq9sec2}), (\ref{eq10sec2}),
we see that $T$ maps the bounded set $\mathcal{M}$ into a compact subset of $%
\mathcal{M}$.

In order to apply the Schauder fixed point theorem we need to prove that the
operator $T$ is continuous. Let $\widetilde{\omega}_{n}, \widetilde{\omega}%
\in {\mathcal{M}}$ and $\widetilde{\omega}_{n}\to\widetilde{\omega}$ in $%
C(0,T;\; L_{2}(\Omega))$. From (\ref{eq66sec2}) and (\ref{ell-2}) of Lemma \ref{teo2sec2}, 
it follows that for any $p>n$
\begin{eqnarray*}
||\nabla h_{n}-\nabla h||_{C^\alpha(\overline\Omega)}&\leqslant& C||[%
\widetilde{\omega}_n]_R-[\widetilde{\omega}]_R||_{L_p(\Omega)}  \notag \\
\leqslant C\cdot(2R)^{\frac{p-2}{p}}||[\widetilde{\omega}_n]_ R-[\widetilde{%
\omega}]_R||^{2/p}_{L_2(\Omega)} &\leqslant& C(R)||\widetilde{\omega}_n-%
\widetilde{\omega}||^{2/p}_{L_2(\Omega)} \overset{n\to\infty}{\longrightarrow%
} 0,
\end{eqnarray*}
or
\begin{eqnarray}
||\nabla h_n-\nabla h||_{C(0,T;\;C^\alpha(\overline\Omega))}\to 0,
\label{eq11sec2}
\end{eqnarray}
where $h_{n}$ and $h $ are the solutions of (\ref{eq4sec2}) with $\omega$
replaced by $\widetilde{\omega}_{n}$ and $\widetilde{\omega}$, respectively.
Next we consider $\omega_n=T [\widetilde{\omega}_n] \mbox{ and } \omega=T [%
\widetilde{\omega}_n]. $ Using \eqref{eq11sec2} and the results of \cite%
{LadySolonUral68} we conclude that
\begin{eqnarray*}
\max_{t\in[0,T]}||\omega_n(\cdot,t)-\omega(\cdot,t)||_{L_2(\Omega)}+
||\nabla(\omega_n-\omega)||_{L_2(0,T;\;L_2(\Omega))}\to 0 .
\end{eqnarray*}
Therefore we conclude that the sequence $\omega_{n}$ itself converges to $%
\omega$ and the continuity of the operator $T$ is proved.

Hence there exists a fixed point $\omega$, such that $\omega=T[\omega].$ We
have shown the following Lemma.

\begin{lemma}
\label{lem3sec2} There exists at least one solution $\{\omega, h\}$ of the
systems (\ref{eq3sec2})-(\ref{eq4sec2}), such that for some $\alpha \in
(0,1) $
\begin{eqnarray*}
\omega \in &W^{1,0}_{2} (\Omega_{T})\cap C(0,T;\; L_{2} (\Omega)),\;\;\;\;\;
h \in &C(0,T;\; C^{1+\alpha } (\overline \Omega)).
\end{eqnarray*}
\end{lemma}

By the theory of parabolic and elliptic equations the constructed functions $%
\omega, h$ have a better regularity. In fact we deduce the following result.

\begin{teo}
\label{teo4sec2} For fixed ${\varepsilon}, R>0$, there exists an unique pair
of functions
\begin{eqnarray}
\omega\in W^{2,1}_{2}(\Omega_{T})\cap
C^{\alpha,\alpha/2}(\Omega_{T}),\,\,\,\,\,\,\,\, h\in
C^{2+\alpha,\alpha/2}(\Omega_{T})  \label{eq20sec2}
\end{eqnarray}
for some $\alpha \in (0,1)$, which is the solution of \textbf{Problem} $%
\mathbf{P_{{\varepsilon},R}}$.
\end{teo}

\noindent\textbf{Proof.} Because of (\ref{eq122sec2}), the function $%
\breve\omega^{\varepsilon}$ is bounded in $\overline\Omega_{T}$ by the
constant $\aleph$. Applying Theorem 7.1, p.181, \cite{LadySolonUral68}, we
have
\begin{eqnarray}
||\omega||_{L_{\infty}(\Omega_{T})} \leqslant C(R, {\varepsilon} )\,\aleph,
\label{eq13sec2}
\end{eqnarray}
where the constant $C(R, {\varepsilon} )$ depends on $R, {\varepsilon}$.
Hence $\omega$ is the solution of the equation
\begin{eqnarray*}
\omega_{t} - {\varepsilon}\Delta\omega={\mathcal{F}}(\vec x,t),\quad(\vec
x,t)\in \Omega_{T}  \label{eq15sec2}
\end{eqnarray*}
with ${\mathcal{F}} = - \nabla\omega\,\vec v - \omega([\omega]_{R}-h)\in
L_{2} (\Omega_{T})$ by (\ref{eq14sec2}), (\ref{eq9sec2}), (\ref{eq13sec2}).
Using Theorem 6.1, p.178, \cite{LadySolonUral68}, we deduce $\;\;\omega(\vec
x,t)\in W^{2,1}_{2} (\Omega_{T}) \;\; \label{eq16sec2} $ and also by Theorem
10.1, p.204, \cite{LadySolonUral68}, we have $\;\;\omega\in C^{\alpha,
\alpha/2}(\Omega_{T}) \;\,\,\,\,\,\, \mbox{ for some  } \alpha \in
(0,1).\;\; \label{eq18sec2}$ Moreover by the theory of elliptic equations
\cite{LadyUral68}, we conclude that $\;\; h\in C^{2+\alpha,
\alpha/2}(\Omega_{T}). \label{eq19sec2}\;\;$

The uniqueness of the solution $\{\omega, h\}$ for \textbf{Problem} $\mathbf{%
P_{{\varepsilon},R}}$ follows in the usual way. Let $\omega
_{i},h_{i},i=1,2,\;$ be different solutions of \textbf{Problem} $\mathbf{P_{{%
\varepsilon},R}}$\ and $\omega=\omega_{1}-\omega_{2},h=h_{1}-h_{2}.\;$ Then
the pair $\{\omega,h\}$ satisfies
\begin{equation*}
\left\{
\begin{array}{cc}
\omega_{t}-{\varepsilon}\Delta\omega=div(\nabla h_{1}\omega+\omega_{2}\nabla
h), \;\;\;\;\;\;\;\;\;\;\;\;\;(\vec{x},t) \in\Omega_{T};   \\
\omega(\vec{x},t)=0, \;\; (\vec{x},t) \in\Gamma_{T};\;\;\;\;\;\;\;\; \omega(%
\vec{x},0)=0, \;\;\vec{x} \in\Omega; &
\end{array}
\right.
\end{equation*}%
\begin{equation*}
\left\{
\begin{array}{cc}
-\Delta h+h=[\omega_{1}]_{R}-[\omega_{2}]_{R}, & (\vec{x},t) \in\Omega_{T},
\\
&  \\
-\nabla h\,\vec{n}=0, & (\vec{x},t) \in\Gamma_{T}.%
\end{array}
\right.
\end{equation*}
Multiplying the first equation by $\omega$, the second one by $h$ and
integrating them over $\Omega$, we obtain the following relations
\begin{eqnarray*}
\frac{1}{2}\frac{d}{dt}\int_{\Omega}\omega^{2}\; d\vec{x}+{\varepsilon}
\int_{\Omega}|\nabla\omega|^{2}\; d\vec{x}&=&-\int_{\Omega}\left( \nabla
h_{1}\omega+\omega_{2}\nabla h\right) \nabla\omega\; d\vec{x}, \\
\int_{\Omega}\left( |\nabla h|^{2}+|h|^{2}\right) \; d\vec{x}&=&\int_{\Omega
}\left( [\omega_{1}]_{R}-[\omega_{2}]_{R}\right) h \;d\vec{x}.
\end{eqnarray*}
Applying the property of \textit{\ the cut-off function}, that $\left|
\lbrack\omega_{1}]_{R}-[\omega_{2}]_{R}\right| \leqslant\left| \omega
_{1}-\omega_{2}\right| , $ and Cauchy's inequality, we come to the
inequality
\begin{eqnarray*}
\frac{d}{dt}\int_{\Omega}\omega^{2}\; d\vec{x}+{\varepsilon}%
\int_{\Omega}|\nabla\omega|^{2}\; d\vec{x} \leqslant
C\int_{\Omega}\omega^{2}\; d\vec{x} \,\,\,\,\,\,\mbox{ and }
\int_{\Omega}\omega^{2}(\vec{x},0)\; d\vec{x} =0
\end{eqnarray*}
with some constant $C=C({\varepsilon},\left\| \nabla h_{1}\right\|
_{L_{\infty}(\Omega_T) },\left\| \omega_{2}\right\| _{L_{\infty
}(\Omega_T)}) $. Application of the standard Gronwall inequality completes
the proof of the uniqueness of solution. $\;\;\;\;\;\;\;\;\blacksquare$

\section{A priori estimates independent of ${\protect\varepsilon}, R$}

In this section we derive a priori estimates of the solution $\{ \omega_{{%
\varepsilon},R}, \,\, h_{{\varepsilon},R} \}$ for \textbf{Problem} $\mathbf{%
P_{{\varepsilon},R}}$, which do not depend on ${\varepsilon}, R$. \label%
{sec4}

\subsection{Maximum-minimum principle for $\protect\omega$ and $h$}

\label{sec4.1} \setcounter{equation}{0}

In this section we show that the solution $\{ \omega_{{\varepsilon},R}, \,\,
h_{{\varepsilon},R} \}$ is bounded in $L_{\infty}(\Omega_{T})$,
independently of ${\varepsilon}, R$. Throughout the section, for simplicity
of presentation, we continue to suppress the dependence of $\omega_{{%
\varepsilon},R}$ and $h_{{\varepsilon},R}$ on ${\varepsilon},R$ and write $%
\omega$ and $h$. All constants $C$ in this section do not depend on ${%
\varepsilon}, R$.

The proof of boundedness of the functions $\omega, h$ is divided on few
lemmas. First let us show the positivity of $\omega$.

\begin{lemma}
\label{lemma1sec4} For all $(\vec x,t)\in\Omega_{T}$
\begin{eqnarray}
\omega(\vec x,t)\geqslant0.  \label{eq1sec3}
\end{eqnarray}
\end{lemma}

\noindent\textbf{Proof.} By (\ref{eq20sec2}) we have
\begin{eqnarray}
\sup_{\vec x\in\Omega}\Big|\mbox{div}\,\vec v(\vec x,t)\Big|\leqslant
\lambda(t)=\max_{\vec x\in\Omega}\Big|\lbrack \omega(\vec x,t)]_{R}-h(\vec
x,t)\Big |\in C(0,T).  \label{eq2sec3}
\end{eqnarray}
Let us denote by $\omega_{-} = \min(\omega,0)$. Then using (\ref{eq122sec2}%
), the first equation of (\ref{eq3sec2}) and the boundary condition $%
\omega_{-} ( \vec x , t )= 0 , \,\,\, \forall(\vec x,t)\in\Gamma_{T}$, it is
easy to verify that the function $\omega_{-}$ satisfies the inequality
\begin{eqnarray*}
\frac{d}{dt}\int_{\Omega}\omega^{2}_{-} \; d\vec{x} + {\varepsilon}%
\int_{\Omega}\Big|\nabla\omega_{-}\Big|^{2} \;d\vec{x}=-\int_{\Omega}%
\mbox{div}\,\vec v\, \omega^{2}_{-} \; d\vec{x}\leqslant\lambda(t)\,
\int_{\Omega}\omega^{2}_{-} \;d\vec{x}.
\end{eqnarray*}
Since $\omega_{-}(\vec x,0)=0,\, \forall\vec x\in\Omega$, from Gronwall
inequality it follows \linebreak $\omega_{-}=0 $ for all $%
(x,t)\in\Omega_{T}, $ which implies (\ref{eq1sec3}). $\;\;\;\;\;\;\;\;%
\blacksquare$\newline

Now we show that $\omega$ is bounded in the space $L_{1} (\Omega)$.

\begin{lemma}
\label{lemma2sec4} There exists a constant $\Upsilon_0, $ independent of $R,
{\varepsilon}$, such that
\begin{eqnarray}
||\omega(\cdot,t)||_{L_{1}(\Omega)}\leqslant \Upsilon_0 ,\quad\forall t\in[%
0,T].  \label{eq3sec3}
\end{eqnarray}
\end{lemma}

\noindent\textbf{Proof.} The function $z=\omega-\breve\omega^{\varepsilon}$
satisfies the problem
\begin{eqnarray}
\left\{
\begin{array}{ll}
z_{t} + \mbox{div}\, (z\vec v)={\varepsilon}\Delta z+F,
\;\;\;\;\;\;\;\;\;\;\; (\vec x,t)\in\Omega_{T}; &  \\
&  \\
z|_{\Gamma_{T}}=0, \;(\vec x,t)\in\Gamma_{T}; \;\;\;\;\;\;z|_{t=0}=0, \;
\vec x\in\Omega &
\end{array}
\right.  \label{eq6sec3}
\end{eqnarray}
with
\begin{eqnarray}
F=({\varepsilon}-1)\partial_{t}{\breve\omega^{\varepsilon}}-\mbox
{div}\,(\vec v \,{\breve\omega^{\varepsilon}}) .  \label{FF}
\end{eqnarray}
By (\ref{eq122sec2}), (\ref{eq1122sec2}) and (\ref{eq4sec2}), we have
\begin{eqnarray}
||\mbox{div}\, (\vec v \,{\breve\omega^{\varepsilon}})(\cdot ,t)||_{L_{1}
(\Omega)} & \leqslant&||{\breve\omega^{\varepsilon}}||_{L_{\infty}(%
\Omega_{T})}\,||[\omega]_{R}-h||_{L_{1}(\Omega)}  \notag \\
+||\nabla{\breve\omega^{\varepsilon} }(\cdot,t)||_{L_{\infty}(\Omega)}
\,||\vec v||_{L_{1} (\Omega)} &\leqslant&
\lambda(t)(||[\omega]_{R}-h||_{L_{1}(\Omega)}+ ||\vec v||_{L_{1}
(\Omega)})\;\;\;  \label{eq7sec3}
\end{eqnarray}
with $\lambda(t)\in L_{1}(0,T)$. Using (\ref{eq66sec2}) and the formulae (\ref{ell-4}), (\ref{ell-5}) 
of Lemma \ref{teo2sec2}, we derive that for
a.e. $t\in[0,T]$,
\begin{eqnarray*}
\left.
\begin{array}{l}
||h(\cdot,t)||_{L_{1}(\Omega)} \\
\\
||\vec v (\cdot,t) ||_{L_{1}(\Omega)} \\
\end{array}
\right\} \leqslant C(||\omega(\cdot,t) ||_{L_{1}(\Omega)}+||a^{{\varepsilon}%
}(\cdot,t) ||_{C^{\eta }(\Gamma)})
\end{eqnarray*}
\begin{eqnarray}
\leqslant C||z (\cdot,t) ||_{L_{1}(\Omega)} + C.\,\,\,\,  \label{eq8sec3}
\end{eqnarray}
Therefore, accordingly $\partial_{t}{\breve\omega^{\varepsilon}}\in L_{1}
(\Omega_{T})$, we deduce that
\begin{eqnarray}
||F(\cdot, t)||_{L_{1}(\Omega)}\leqslant\lambda(t)\left( ||z (\cdot,t)
||_{L_{1}(\Omega)} + 1\right)\,\,\, \mbox{ with } \,\,\,\lambda(t)\in
L_{1}(0,T).  \label{eq9sec3}
\end{eqnarray}

Multiplying the equation of (\ref{eq6sec3}) by $sgn_{\delta}z:=\frac {z}{%
\sqrt{z^{2}+\delta}}$ with some $\delta \in (0,1)$, we obtain that
\begin{eqnarray}
\partial_t ( (z^{2}+\delta)^{1/2} ) + \text{div} (\vec
v\,(z^{2}+\delta)^{1/2})&-&\frac{\delta }{(z^{2}+\delta)^{1/2}}\,\text{div}%
\vec v  \notag \\
={\varepsilon}\,\text{div}\left( \nabla z \,\frac {z }{(z ^{2}+\delta)^{1/2}}%
\right) &-&{\varepsilon} \delta\frac{|\nabla z |^{2}}{(z ^{2}+\delta)^{3/2}}%
+F \, sgn_{\delta }z.\;\;\;\;\;\;  \label{eqAB7sec3}
\end{eqnarray}
Taking into account that $z=0 $ on $\Gamma_{T}$ and integrating this
equality over $\Omega$, we have
\begin{eqnarray*}
\frac{d}{dt}\Big(||\sqrt{z^{2}+\delta}||_{L_{1}(\Omega)}\Big)+I_{\delta
}(t)\leqslant\lambda(t)\,\left( ||z||_{L_{1}(\Omega)}+1\right)
\mbox{ for
a.e. } t\in[0,T],
\end{eqnarray*}
with $I_{\delta}(t)=\int_{\Omega} \Bigl[ \mbox{div}(\vec v\,\sqrt {%
z^{2}+\delta})-\frac{\delta}{(z^{2}+\delta)^{1/2}}\,\mbox{div}\,\vec v %
\Bigl] \; d\vec{x}. \,\, $ Since
\begin{eqnarray*}
|I_{\delta}|\leqslant\delta^{1/2}\,\{ \int_{\Omega}(|\omega|+|h|)\; d\vec{x}
+ | \int_{\Gamma}a^{\varepsilon} \; d\vec{x} | \} \mathop{\longrightarrow}%
\limits_{{\delta}\to 0} 0,
\end{eqnarray*}
we deduce the inequality
\begin{eqnarray*}
\frac{d}{dt}||z||_{L_{1}(\Omega)}\leqslant\lambda(t)\,(
||z||_{L_{1}(\Omega)}+1 ), \mbox{ a.e. } t\in[0,T].
\end{eqnarray*}
Hence applying Gronwall inequality, we obtain $||z(\cdot,t)
||_{L_{1}(\Omega)}\leqslant C \mbox{ for } t\in[0,T], $ that immediately
implies the estimate (\ref{eq3sec3}). $\;\;\;\;\;\;\;\;\blacksquare$\newline

Now we prove an auxiliary Lemma.

\begin{lemma}
\label{lemma3sec4} There exists a positive constant $R_{\ast }$, such that  
\begin{equation}
\max_{\overline{\Omega }_{T}}\omega (\vec{x},t)\leqslant max\{\,\max_{%
\overline{\Omega }_{T}}h(\vec{x},t),\;\aleph \,\} \label{eq10'sec3}
\end{equation}
for any fixed  
$ R>R_{\ast }.$
The constant $R_{\ast }$ depend on $\Upsilon _{0},\Gamma ,\Omega $ and $%
\aleph ,$ being independent of ${\varepsilon .}$
\end{lemma}

\noindent \textbf{Proof.} A) If $\omega $ attains its maximum on the
parabolic boundary $\Gamma_T$, then (\ref{eq10'sec3}) is obvious.

B) If $\omega $ attains a positive maximum $\overline{\omega }=\omega (\vec{x%
}_{0},t_{0})$ inside $\Omega _{T}$, then we come to the relation%
\begin{equation}
0\leqslant \omega _{t}(\vec{x}_{0},t_{0})-\varepsilon \vartriangle \omega (%
\vec{x}_{0},t_{0})= \overline{\omega }(h(\vec{x}_{0},t_{0})-[\overline{%
\omega }]_{R})  \label{add-1}
\end{equation}
and hence 
\[
[\overline{\omega }]_{R}=\min \{R,\overline{\omega }\}\leqslant \max_{%
\overline{\Omega }_{T}}h(\vec{x},t) .
\]
Let us consider two possibilities:

B.1) If $\overline{\omega }\leqslant R\ $, we have $\ [\overline{\omega }%
]_{R}=\overline{\omega }\leqslant \max_{\overline{\Omega }_{T}}h(\vec{x},t),$
therefore (\ref{eq10'sec3})\ is proved.

B.2) If $\overline{\omega }>R$, we have 
\begin{equation}
R\leqslant \max_{\overline{\Omega }_{T}} h(\vec{x},t).  \label{add-11}
\end{equation}
Let us prove that this case is impossible for large enough $R$. Using the representation 
(\ref{eq66sec2}) and (\ref{ell-1}), 
(\ref{ell-4}), (\ref{ell-5}) of
 Lemma \ref{teo2sec2}, we have 
\begin{eqnarray*}
\max_{\overline{\Omega }}\,|h| &\leqslant &\max_{\overline{\Omega }%
}\,|K_{1}\ast ({[\omega ]_{R}}^{1/q}{[\omega ]_{R}}^{1-1/q})|\,\,+\,\,\max_{%
\overline{\Omega }}\,|h_{2}| \\
&\leqslant &C(q,\Omega )\,\,||{\omega }||_{L_{1}(\Omega )}\,\,R^{1-1/q}\,\,+\,\,C(p,\Gamma )||a^{\varepsilon
}||_{C^{\eta}(\Gamma )} \\
&\leqslant &C(q,\Omega )\,\,\Upsilon _{0}\,\,R^{1-1/q}\,\,+\,\,C(p,\Gamma )||a^{\varepsilon }||_{C^{\eta}(\Gamma )}
\end{eqnarray*}%
with $n<q<\infty ,\;n-1<p<\infty .$ Last inequality implies that 
\[
\max_{\overline{\Omega }_{T}}\,|h|\leqslant C_{\ast }(\max_{\overline{\Omega 
}_{T}}\,R^{1-1/q}+1)
\]%
with some constant $C_{\ast }=C(q,\Omega ,p,\Gamma ,\Upsilon
_{0},||a||_{L_{\infty }(0;T;\,\,C^{\eta}(\Gamma ))}),$ which does not depend on $%
\varepsilon ,R.$ Since $\;\;0<1-1/q<1, \;\;$ there exists a constant $R_{\ast },$ depending on $%
C_{\ast },$ such that $\;\;C_{\ast }(\,R^{1-1/q}+1)%
\leqslant \frac{R}{2}\;\;$ for any $\;\;R>R_{\ast }.\;\;$ Therefore we conclude%
\[
\max_{\overline{\Omega }_{T}}\,|h|\leqslant \frac{R}{2}\;\;\;\;\mbox{ for any  }\;\;\;\;%
R>R_{\ast }.
\]%
But this inequality contradicts with (\ref{add-11}).$\;\;\;\;\blacksquare $

\begin{remark}
Let us note  that if consider the system (\ref{eq3sec2})-(\ref{eq4sec2}) with the right side in the equation of 
(\ref{eq4sec2}) equals to $\omega$, instead of $[\omega ]_{R}$, then the desired assertion of Lemma \ref{lemma3sec4} 
 follows immediately
from an analog of the relation (\ref{add-1})
\begin{equation}
0\leqslant \omega _{t}(\vec{x}_{0},t_{0}) - \varepsilon \vartriangle \omega (\vec{x}_{0},t_{0})=
\overline{\omega }(h(\vec{x}_{0},t_{0})-\overline{\omega }).
\end{equation} \notag
\end{remark}

Now we are able to show the boundedness of $\{ \omega, h \} $ in $%
L_{\infty}(\Omega_{T} ) .$

\begin{lemma}
\label{lem4sec3} There exist constants $\Upsilon_{1}, \Upsilon_{2} $
depending on the data $\mathbf{\breve\omega}$, $a$, but independent of $R, {%
\varepsilon}$, such that
\begin{eqnarray}
||\omega||_{L_{\infty}(\Omega_{T})}&\leqslant&\Upsilon_{1},
\label{eq13+sec3} \\
||h||_{L_{\infty}(\Omega_{T})}, \;||\nabla
h||_{L_{\infty}(\Omega_{T})}&\leqslant&\Upsilon_{2} .  \label{eq99sec3}
\end{eqnarray}
\end{lemma}

\noindent\textbf{Proof. } Using the representation (\ref{eq66sec2}), Lemma \ref{lemma2sec4} and (\ref{ell-1}), 
(\ref{ell-4}), (\ref{ell-5}) of
 Lemma \ref{teo2sec2},  we have
\begin{eqnarray*}
\max_{\overline\Omega} \,|h|&\leqslant& \max_{\overline\Omega} \,|K_1\ast({%
[\omega]_{R}}^{1/q}{[\omega]_{R}}^{1-1/q})| \,\,+\,\,\max_{\overline\Omega}
\,|h_2|  \label{eq13''''sec3} \\
&\leqslant& C(q,\Omega)\,\,\Upsilon_0 \,\, \max_{\overline\Omega}
\;\omega^{1-1/q} \,\,+\,\,C(p,\Gamma)||a^{\varepsilon}||_{C^{\eta}(\Gamma)}
\end{eqnarray*}
with $n<q<\infty,\;n-1<p<\infty.$ Last inequality implies that
\begin{equation}  \label{+}
\max_{\overline\Omega_T}\,|h|\leqslant
C(\max_{\overline\Omega_T}\,\omega^{1-1/q}
+||a||_{L_{\infty}(0;T; \,\,C^{\eta}(\Gamma))})
\end{equation}
with some constant $C=C(q,\Omega,p,\Gamma,\Upsilon_0 ),$ which does not
depend on $\varepsilon, R.$ Joining (\ref{eq10'sec3}) and (\ref{+}), we
obtain the inequality
\begin{equation}  \label{++}
\max_{\overline\Omega_T}\,\omega \leqslant C \max\{
\,\max_{\overline\Omega_T}\,\omega^{1-1/q}
+||a||_{L_{\infty}(0;T;\,\,C^{\eta}(\Gamma))},\,\,\aleph\, \},
\end{equation}
which leads to (\ref{eq13+sec3}), because $\,\,0<1-1/q<1.$ The estimates (\ref{eq99sec3})
 follow immediately from (\ref{ell-2}), (\ref{ell-5}) of Lemma \ref{teo2sec2}. This completes
the proof. $\;\;\;\;\;\;\;\;\blacksquare$\newline

Choosing $R:=\Upsilon_{1} $ in (\ref{eq4sec2}), we see that the cut-off
function $[\cdot]_{R}$ in (\ref{eq4sec2}) can be omitted. In the following
we consider the solution of (\ref{eq3sec2}), (\ref{eq4sec2}) as $%
\omega_{\varepsilon},\,h_{\varepsilon}$, depending only on ${\varepsilon}$
(not $R$).

\subsection{Estimates of derivatives}

\label{sec4.2}

From this section we shall write $\omega _{\varepsilon}$, $h_{\varepsilon} $
and $\vec{v}_{\varepsilon}=-\bigtriangledown h_{\varepsilon}$. Let us show
the following lemma.

\begin{lemma}
\label{lemma4sec4} There exist constants $C$ independent of ${\varepsilon}$,
such that
\begin{eqnarray}
||\sqrt{{\varepsilon}}\nabla\omega_{\varepsilon}||_{L_{2}
(\Omega_{T})}&\leqslant& C,  \label{eq23sec3} \\
||\partial_{t}(\omega_{\varepsilon})||_{L_{2} (0,T; \,\,H^{-1}(\Omega
))}&\leqslant& C,  \label{eq30sec3} \\
||h_{\varepsilon}||_{L_\infty (0,T; \; C^{1+\eta } (\overline{\Omega}) )}
&\leqslant& C,  \label{eq255sec3} \\
||h_{\varepsilon}(\cdot,t)||_{L_\infty (0,T; \; W^{2}_{p} (\overline{%
\Omega^{\prime }}) ) } &\leqslant& C,  \label{eq955sec3}
\end{eqnarray}
for any $p\in[1,\infty)$ and for any $%
\Omega^{\prime }$, such that $\overline{\Omega^{\prime }}\subset \Omega$.
\end{lemma}

\noindent\textbf{Proof.} According to Lemma \ref{lem4sec3}, we have, that
\begin{equation}
||h_{\varepsilon}||_{L_{\infty}(\Omega_{T})}, \,\,||\nabla h_{\varepsilon}
||_{L_{\infty}(\Omega_{T})} ,
\,\,||\omega_{\varepsilon}||_{L_{\infty}(\Omega_{T})} \leqslant\Upsilon_{3},
\label{eq26sec3}
\end{equation}
where $\Upsilon_{3}$ is independent of ${\varepsilon}$. Multiplying the
equation of \eqref{eq6sec3} by $z_{\varepsilon}=\omega_{{\varepsilon}%
}-\breve\omega^{\varepsilon}$ and using (\ref{eq122sec2}), (\ref{eq26sec3}),
$div\,\vec{v}_{\varepsilon}=\omega_{\varepsilon}-h_{\varepsilon}$, we easily
get
\begin{eqnarray*}
\frac{1}{2}\frac{d}{dt}||z_{\varepsilon}||^{2}_{L_{2} (\Omega)}+ {\varepsilon%
}||\nabla z_{\varepsilon}||_{L_{2} (\Omega)} & \\
\leqslant -\int_{\Omega}\mbox{div}\,\vec v_{\varepsilon}\frac{%
z_{\varepsilon}^{2}}{2}\; d\vec{x} &+\int_{\Omega}F_{\varepsilon} \;
z_{\varepsilon} \; d\vec{x} \leqslant C+ C||F_{{\varepsilon} }||_{L_{1}
(\Omega)}.
\end{eqnarray*}
Integrating on $(0,T)$ and taking into account (\ref{eq03sec1}), (\ref%
{eq1122sec2}), we have
\begin{eqnarray}
||\sqrt{{\varepsilon}}\nabla z_{\varepsilon}||_{L_{2}(\Omega_{T})}\leqslant
C,  \label{44.5}
\end{eqnarray}
that implies (\ref{eq23sec3}).

Let us choose an arbitrary function $\phi \in H^{1}(\Omega _{T})$, such that
$\phi (\vec{x},T)=\phi(\vec{x},0)=0$ for all $\vec{x}\in \Omega $ and $\phi (%
\vec{x},t)=0, \; (\vec{x},t)\in \Gamma_T$. Multiplying the first equation in
(\ref{eq3sec2}) by $\phi$ and integrating it over $\Omega_T$, from (\ref%
{eq26sec3})-(\ref{44.5}) we obtain
\begin{equation*}
|\int\limits_{\Omega_T }\omega _{{\varepsilon} }\,\phi _{t}\; d\vec{x}dt
\,|\leqslant C\,\Vert \phi \Vert _{L_{2}(0,T;\; H_0^{1}(\Omega )) },
\end{equation*}
that gives (\ref{eq30sec3}).

Taking into account
the representation  (\ref{eq66sec2}), we have
\begin{equation}
h_{{\varepsilon}}(\vec{x},t)=(K_{1} \ast \omega_{\varepsilon}) (\vec{x}%
,t)+(K_{2} \ast a^{\varepsilon}) (\vec{x},t).  \label{present}
\end{equation}
Hence by
 \eqref{ell-21}, \eqref{ell-5} of Lemma \ref{teo2sec2}
 we have for all $t\in(0,T)$
\begin{equation*}  \label{time-5}
||h_{\varepsilon}(\cdot,t)||_{C^{1+
\eta } (\overline{\Omega})}
\leqslant C\left(||\omega_{\varepsilon}(\cdot,t)||_{L_{\infty}
(\Omega)}+||a^{\varepsilon}(\cdot,t)||_{C^{\eta } (\Gamma)}\right)
\end{equation*}
and by \eqref{ell-3}, \eqref{ell-6} of Lemma \ref{teo2sec2}
\begin{equation*}  \label{time-6}
||h_{\varepsilon}(\cdot,t)||_{W_{p}^{2} (\overline{\Omega^{\prime }})}
\leqslant C\left(||\omega_{\varepsilon}(\cdot,t)||_{L_{\infty}
(\Omega)}+||a^{{\varepsilon}}(\cdot,t)||_{L_{1} (\Gamma)}\right) ,\;\,\,
\end{equation*}
for all $p\in [1,\infty)$ and for any $\Omega^{\prime }$ such that $%
\overline{\Omega^{\prime }} \subset \Omega $. With the help of \eqref{eq00sec1}, \eqref{ABsec2} and \eqref{eq26sec3}, we
derive the assertions (\ref{eq255sec3}), (\ref{eq955sec3}). $%
\;\;\;\;\;\;\;\;\blacksquare$

\section{Limit transition}

\label{sec5}

In this section we prove Theorem 2.

From (\ref{eq955sec3}), (\ref{eq26sec3}) we conclude that there exists a
subsequence of $\{\omega_{\varepsilon},h_{\varepsilon} \}$ such that
\begin{eqnarray}
h_{\varepsilon} \rightharpoonup h &\mbox{ weakly}-\ast& \mbox{ in } L_\infty
(\Omega_T ) \cap L_\infty (0,T;\, W^2_p(\overline{\Omega^{\prime }} )),
\notag \\
\nabla h_{\varepsilon} \rightharpoonup \nabla h &\mbox{ weakly}-\ast&
\mbox{
in } L_\infty (\Omega_T ) ,  \label{eq1032sec3}
\end{eqnarray}
for any $p \in(1,\infty)$ and for any $\Omega^{\prime }$, such that $%
\overline{\Omega}^{\prime }\subset \Omega .$ By results obtained in \cite%
{aubin}, \cite{simon}, the compact embedding of $L_2(\Omega)$ in $%
H^{-1}(\Omega)$ and (\ref{eq30sec3}) imply
\begin{equation}
\omega _{{\varepsilon} }\to \omega \quad \mbox{strongly in}\quad L_2 (0,T ;
\; H^{-1}(\Omega ) ).  \label{4.9}
\end{equation}
In view of (\ref{eq23sec3}), (\ref{eq26sec3}), we have
\begin{eqnarray}
\omega_{\varepsilon}\rightharpoonup\omega &\mbox{ weakly}-\ast& \mbox{ in }
L_{\infty}(\Omega_{T}),  \notag \\
{\varepsilon} \nabla \omega_{\varepsilon} \rightharpoonup 0 &
\mbox{ weakly
in } L_{2}(\Omega_{T}).  \label{eq31sec3}
\end{eqnarray}

With the help of (\ref{ABsec2}), (\ref{eq1032sec3}), (\ref{eq31sec3}), the
limit transition in the representation (\ref{present}) on ${\varepsilon}\to 0
$ yields
\begin{equation}
h(\vec{x},t)=(K_{1} \ast \omega) (\vec{x},t)+(K_{2} \ast a) (\vec{x},t)
\,\,\, \mbox{ for a.e. } (\vec x ,t)\in \Omega_T .  \label{Apresent}
\end{equation}
From the elliptic theory \cite{LadyUral68}, \cite{Vlad00} and (\ref{ell-21}), (\ref{ell-5}) of Lemma \ref%
{teo2sec2}, the function $h$ satisfies the equation (\ref{eq8}) with the
boundary conditions (\ref{eq10}), such that
\begin{equation}
h\in L_\infty (0,T;\, C^{1+\eta } (\overline{\Omega} )). \label{Bpresent}
\end{equation}

\begin{definition}
\label{def2sec5} Let the distance between any given point $\vec x \in\mathbb{%
R}^{n}$ and any subset $A \subseteq\mathbb{R}^{n} $ be defined by $d(\vec x,
A):=inf_{\vec y \in A } |\vec x- \vec y| $. Let $d=d( \vec x )$ be the
distance function on $\Gamma$, defined by
\begin{eqnarray*}
d(\vec x):=d(\vec x,\mathbb{R}^{n}\backslash\Omega)-d(\vec x,\Omega) \,\,\, %
\mbox{ for any } \vec x \in\mathbb{R}^{n} .
\end{eqnarray*}
Let $A$ be an arbitrary sub domain of $\Omega $. We introduce also the
distance between $A$ and the boundary $\Gamma$ by $d(A , \Gamma ):=inf_{\vec
y \in A } d(\vec y) . $
\end{definition}

The set of all points of $\overline{\Omega}$, whose distance to $\Gamma$ (
to $\Gamma^-$ and to $\Gamma^+$) is less than $\sigma$, is denoted by $%
U_{\sigma}(\Gamma)$ ( respectively, $U_{\sigma}(\Gamma^-)$ and $%
U_{\sigma}(\Gamma^+)$). Since $\Gamma\in C^{2+\gamma}$, the function $d=d(
\vec x )$ belongs to $C^{2}$ in a neighbourhood $U_{\sigma_{0}}(\Gamma)$ of $%
\Gamma$ for some ${\sigma_{0}}>0$. For $0<2 {\sigma}<{\sigma_{0}}$, we
introduce the approximation of the unit function for all $\vec x \in
\overline \Omega $ by
\begin{eqnarray}
{\mathbf{1}_{\sigma}} (\vec x):=\left\{
\begin{array}{ll}
1, & \mbox{ if } \vec x\in \Omega \backslash U_{2\sigma}(\Gamma), \\
&  \\
\frac{d-\sigma}{\sigma}, & \mbox{ if } \sigma <d(\vec x) < 2\sigma, \\
&  \\
0, & \mbox{ if } 0 \leqslant d(\vec x) < \sigma .%
\end{array}
\right.  \label{ddeq6sec3}
\end{eqnarray}

Now we show an auxiliary lemma, playing the crucial role in the proof that
the function $\omega$ satisfies the boundary condition (\ref{eq5}) in the
sense of the equality (\ref{eq7}).

\begin{lemma}
For any positive $\psi\in C^{1,1} (\Omega_{T})$, such that $supp
\,(\psi)\subset \Omega_{T}\cup\Gamma^{-}_{T},$ and $\psi(\vec x,T)=0, \,\,
\vec x\in\Omega$, we have
\begin{eqnarray}
\lim_{{\sigma}\to0} \,\,\left(\, \overline{ \lim_{{\varepsilon}\to0} }\, \,%
\frac{1}{{\sigma}} \int^{T}_{0}\int_{[\sigma<d<{2\sigma}]}|\omega_{%
\varepsilon}-\breve\omega_{\varepsilon} | \,\,\vec v _{\varepsilon}
\bigtriangledown d \,\,\psi \; d\vec{x} dt \, \right) =0.  \label{eq35sec3}
\end{eqnarray}
\label{lem6sec4}
\end{lemma}

\noindent\textbf{Proof.} Taking into account that $\,\,\,z_{\varepsilon}=%
\omega_{\varepsilon}-\breve\omega_{\varepsilon}=0\,\,\,$ on $\Gamma_{T}$ and
multiplying the equality (\ref{eqAB7sec3}) by an arbitrary non negative
function $\eta\in H^{1}(\Omega_{T})$ with $\eta(\vec x,T)=0, \, \vec x
\in\Omega$, we derive that
\begin{eqnarray*}
-\int\int_{\Omega_{T}}\eta_{t} (z_{\varepsilon}^{2}+\delta)^{1/2} \; d\vec{x}
dt&+&\int_{\Omega} \eta(\vec x,0) (z_{\varepsilon}^{2}(\vec
x,0)+\delta)^{1/2} \; d\vec{x} \\
+\int_{\Gamma_{T}} a^{\varepsilon}\,\sqrt{\delta}\,\eta\; d\vec{x} dt
&-&\int\int_{\Omega_{T}} (z_{{\varepsilon}}^{2}\,+ \delta)^{1/2} \, \vec
v_{\varepsilon}\nabla\eta\; d\vec{x} dt \\
-\int\int_{\Omega_{T}}\frac{\delta }{(z_{\varepsilon}^{2}+ \delta)^{1/2}}%
\,(\omega_{{\varepsilon}}-h_{\varepsilon})\eta \; d\vec{x} dt &=&-{%
\varepsilon}\int\int_{\Omega_{T}}\nabla z_{\varepsilon}\nabla\eta \,\frac{%
z_{\varepsilon}}{(z_{\varepsilon}^{2}+\delta)^{1/2}} \; d\vec{x} dt \\
-{\varepsilon}\delta\int\int_{\Omega_{T}}\frac{|\nabla z_{\varepsilon}|^{2}}{%
(z_{\varepsilon}^{2}+\delta)^{3/2}}\,\eta\; d\vec{x} dt
&+&\int\int_{\Omega_{T}}F_{\varepsilon}\, sgn_{\delta}z_{{\varepsilon}%
}\,\eta\; d\vec{x} dt.
\end{eqnarray*}
Let us denote by $G_{\varepsilon}:=|\partial_t \breve{ \omega^{\varepsilon} }
|+|\bigtriangledown \breve{\omega^{\varepsilon}}|$ and $G:=|\partial_t
\breve\omega|+|\bigtriangledown {\breve\omega}|.$ Since the functions $\{
z_{\varepsilon}, \, \omega_{\varepsilon} , \, h_{\varepsilon} \}$ and $a^{{%
\varepsilon}}$ are uniformly bounded with respect to ${\varepsilon}$ in $%
L_\infty (\Omega_T )$ and in $L_\infty (\Gamma_{T})$, using (\ref{FF}) and (%
\ref{eq23sec3}), we obtain
\begin{eqnarray*}
-\int\int_{\Omega_{T}} ( z_{{\varepsilon}}^{2}+\delta)^{1/2}\; \vec
v_{\varepsilon}\nabla\eta \; d\vec{x} dt&\leqslant&
C\int\int_{\Omega_{T}}(|\eta _{t}|+|\eta|+ |\eta|\,|G_{\varepsilon}|) \; d%
\vec{x} dt \\
+C \int_{\Omega} |\eta(\vec x,0) | \; d\vec{x} &+&C\,\sqrt{\delta}%
\int_{\Gamma_{T}} |\eta| \; d\vec{x} dt+C\sqrt {\varepsilon}
||\nabla\eta||_{ L_2 ( \Omega_T )} .
\end{eqnarray*}
Taking $\delta\to 0$, we get
\begin{eqnarray}
-\int\int_{\Omega_{T}} |z_{{\varepsilon}}| \,\,\vec
v_{\varepsilon}\nabla\eta \; d\vec{x} dt&\leqslant&
C\int\int_{\Omega_{T}}(|\eta_{t}|+|\eta|+|\eta |\,|G_{\varepsilon}|)\; d\vec{%
x} dt  \notag \\
&+&C \int_{\Omega} |\eta(\vec x,0) | \; d\vec{x}+C\sqrt {\varepsilon}
||\nabla\eta||_{ L_2 ( \Omega_T )} .\;\;\;\;\;\;  \label{eq36sec3}
\end{eqnarray}
Let us choose an arbitrary function $\psi\in C^{1,1}(\Omega_{T})$,
satisfying the conditions of the present lemma and take the test function in
(\ref{eq36sec3}) as $\eta_{\sigma}= ( 1-{\mathbf{1}_{\sigma}} )\,\psi $ we
derive
\begin{eqnarray*}
\int \int_{\Omega_{T}}|z_{{\varepsilon}}|\,\psi\,\frac{\vec
v_{\varepsilon}\nabla d}{{\sigma}}\, \chi_{[\sigma<d<{2\sigma}]}\; d\vec{x}
dt & \\
\leqslant C\int\int_{\Omega_{T}} ( |\partial_{t} \eta_{{\sigma}}|&+&
|\eta_{\sigma}|+|\eta_{\sigma}|\,|G_{\varepsilon}|+( 1-{\mathbf{1}_{\sigma}}
)|\nabla \psi| ) \; d\vec{x} dt \\
&+&C \int_{\Omega} |\eta_{\sigma} (\vec x,0) | \; d\vec{x}+C\sqrt {%
\varepsilon} ||\nabla\eta_\sigma ||_{ L_2 ( \Omega_T )} .
\end{eqnarray*}
By (\ref{eq255sec3}) the set $\{\vec v_{\varepsilon} ( \cdot , t)\} $ is
uniformly continuous on $\overline{\Omega}$, independently of ${\varepsilon}
$ and $t\in [0,T]$. Also we have that the function $d(\vec x)\in C^{2}$ in
the neighbourhood $U_{\sigma_{0}}(\Gamma)$ and $\bigtriangledown d =-\vec n$
on $\Gamma$. Hence there exists $\sigma _1<\sigma_{0}$, independent of ${%
\varepsilon}$ and $t\in [0,T],$ such that
\begin{eqnarray}
\bigl\{\psi\,\,\vec v _{\varepsilon} \bigtriangledown d \, \bigl\} \,(\vec
x, t) =\left\{
\begin{array}{ll}
>0, & \mbox{ if } (\vec x, t)\in U_{\sigma_1}(\Gamma^-)\times [0,T]; \\
&  \\
0, & \mbox{ if } (\vec x, t)\in U_{\sigma_1}(\Gamma^+ \cup \Gamma^0)\times [0,T] .%
\end{array}
\right.  \label{B1}
\end{eqnarray}
Therefore for any $2 \sigma \in (0,\sigma_1)$, according (\ref{eq1122sec2}),
we have
\begin{eqnarray*}
0 &\leqslant & \overline{\mathop{\lim} \limits_{{\varepsilon}\to0}} \,\,
\frac 1 \sigma \int_0^T \int_{[\sigma<d<{2\sigma}]} \psi\,\,|z_{{\varepsilon}%
}|\,\vec v_{\varepsilon}\nabla d \; d\vec{x} dt \leqslant C \int_{\Omega}
|\eta_{\sigma} (\vec x,0) |\; d\vec{x} \\
&+& C\int\int_{\Omega_{T}} ( |\partial_{t} \eta_{{\sigma}}|+
|\eta_{\sigma}|+|\eta_{\sigma}|\,|G|+( 1-{\mathbf{1}_{\sigma}} )|\nabla
\psi| ) \; d\vec{x} dt \; d\vec{x}.
\end{eqnarray*}
This implies the property (\ref{eq35sec3}), since $\eta_{\sigma}=( 1-{%
\mathbf{1}_{\sigma}} )\psi $, $\partial_t \eta_{{\sigma}}=( 1-{\mathbf{1}%
_{\sigma}} )\psi_{t}$ and $( 1-{\mathbf{1}_{\sigma}} )\to0$ in $\Omega_{T}$,
when ${\sigma} \to0$. $\;\;\;\;\;\;\;\;\blacksquare$\newline

Now we are able to prove that the pair $\{\omega,h\}$ satisfies (\ref{eq7}).
For any $2 \sigma\in (0, \sigma_1 )$, we put $\eta_{\sigma}= {\mathbf{1}%
_{\sigma}}\,\psi, $ where $\psi$ is an arbitrary function, satisfying the
conditions of Lemma \ref{lem6sec4}. Clearly, $\eta_{\sigma}\in
H^1(\Omega_{T}), \,\, \eta_{\sigma}(\vec x, T)=0$ for $\vec x \in \Omega $
and $\eta_{\sigma}=0$ on $\Gamma_{T}$. Multiplying the first equation in (%
\ref{eq3sec2}) by $\eta_{\sigma}$ and integrating it over $\Omega_T$, we
derive
\begin{eqnarray}
0= \{ \int_{\Omega_{T}} [ \omega_{\varepsilon}(\psi_{t}&+&\vec
v_{\varepsilon} \nabla \psi) ] {\mathbf{1}_{\sigma}} -{\varepsilon} \nabla
\omega_{\varepsilon} \nabla \eta_\sigma \; d\vec{x} dt + \int_{\Omega}\breve
\omega^{\varepsilon}(\vec x, 0) \, \eta_\sigma(\vec x, 0) \;d\vec{x} \,\}
\notag \\
&+&\frac {1}{{\sigma}}\int_{0}^{T}\int_{[\sigma<d<{2\sigma}]}
\omega_{\varepsilon} \, (\vec v_{\varepsilon} \nabla d) \, \psi\; d\vec{x}
dt = I^{{\varepsilon},\sigma}+J^{{\varepsilon},\sigma}.  \label{ddsec3}
\end{eqnarray}
Using (\ref{eq1032sec3}), (\ref{4.9}) and $\vec{v}_{\varepsilon} \nabla \psi
\, {\mathbf{1}_{\sigma}} \mathop{\rightharpoonup}\limits_{{\varepsilon}\to0}
\vec{v} \nabla \psi \, {\mathbf{1}_{\sigma}} $ weakly in $%
L_{2}(0,T;\;H_0^1(\Omega))$ with $\vec v=-\nabla h$, we have
\begin{equation*}
\lim\limits_{{\varepsilon} \to 0} \int\limits_{\Omega_T } \omega
_{\varepsilon}\, \vec{v}_{\varepsilon} \nabla \psi \, {\mathbf{1}_{\sigma}}
\; d\vec{x}dt =\int\limits_{\Omega_T}\omega \,\vec{v} \nabla \psi \, {%
\mathbf{1}_{\sigma}} \; d\vec{x}dt.
\end{equation*}
Because ${\mathbf{1}_{\sigma}}\mathop{\longrightarrow}\limits_{{\sigma}\to0}
1$ in $\Omega_{T}$ and $\Omega$, we get
\begin{eqnarray*}
\lim_{{\sigma}\to0} \,\,\left(\, \lim_{{\varepsilon}\to0} I^{{\varepsilon}%
,\sigma}\right) =\int_{\Omega_{T}}\omega(\psi_{t}+\vec v\nabla\psi)\; d\vec{x%
} dt+\int_{\Omega}\omega_{0}\psi(\vec x, 0) \;d\vec{x}.
\end{eqnarray*}
Since
\begin{eqnarray*}
J^{{\varepsilon},\sigma} = \Bigl[\frac{1}{{\sigma}}\int^{T}_{0}\int_{[%
\sigma<d<{2\sigma}]}z_{\varepsilon} \,(\vec v_{\varepsilon} \nabla
d)\,\psi\; d\vec{x} dt\Bigl] +& \\
+ \Bigl[\frac{1}{{\sigma}}\int^{T} _{0}\int_{[\sigma<d<{2\sigma}]}&
\breve\omega^{\varepsilon} \,(\vec v_{\varepsilon} \nabla d)\, \psi\;d\vec{x}%
dt \Bigl]=J_{1}^{{\varepsilon},\sigma}+J_{2}^{{\varepsilon},\sigma}.
\end{eqnarray*}
According to Lemma {\ref{lem6sec4}}, we have
\begin{eqnarray*}
\lim_{{\sigma}\to0} \,\,\left(\, \overline{ \lim_{{\varepsilon}\to0} }
|J_{1}^{{\varepsilon},\sigma}|\right)=0 .
\end{eqnarray*}
By (\ref{Bpresent}) the set of functions $\vec v_{\varepsilon}
\bigtriangledown d $ is uniformly continuous on $\overline{\Omega}$ for all $%
t\in [0,T]$, independently of ${\varepsilon}$ and the trace of $\vec
v_{\varepsilon} \bigtriangledown d $ on $\Gamma_T$ satisfies in the usual
sense
\begin{eqnarray*}
\vec v_{\varepsilon} \bigtriangledown d =- a^{\varepsilon} \, \mbox{ for }
(\vec x,t)\in\Gamma_{T}.
\end{eqnarray*}
By the standard theory of traces and (\ref{eq01sec1}), (\ref{eq03sec1}), (%
\ref{eq1122sec2}), the function $\breve\omega_{\varepsilon}$ has a trace on
the boundary $\Gamma^{-}_{T}$, which converges to $b$, when ${\varepsilon}
\to 0$, in the space $L_1 (0,T;\,L_q (\Gamma^- ))$ for any $q\in[1,\infty)$.
Therefore the convergences (\ref{ABsec2}), (\ref{eq1122sec2}) imply
\begin{eqnarray*}
\lim_{{\sigma}\to0} \,\,\, ( \lim_{{\varepsilon} \to0} \,\, J_{2}^{{%
\varepsilon},\sigma})=-\int^{T}_{0}\int_{\Gamma^{-}_{T}}ab\,\psi\; d\vec{x}
dt.
\end{eqnarray*}
Therefore the pair $\{\omega, \, h\}$ satisfies the equation (\ref{eq7}) for
any $\psi$, satisfying the conditions of Lemma \ref{lem6sec4}. In view of
the linearity of (\ref{eq7}) with respect to $\psi$, we see that this
equation is fulfilled also for any $\psi$, satisfying the conditions (\ref%
{psi}) of Definition \ref{def1sec1}.

Let us show that the limit functions $\omega, h, \bigtriangledown h$ admit
additional regularity with respect to the time variable $t$.

\begin{lemma}
\label{lemma44} For any $\Omega^{\prime }$, such that $\overline{\Omega}%
^{\prime }\subset \Omega$, we have
\begin{eqnarray}
\partial_{t} \, \omega &\in & L_{\infty} (0,T; \,\,H^{-1}(\Omega )),
\label{A144} \\
h, \,\, \bigtriangledown h &\in & C^{0 , \, \theta } (\overline{%
\Omega^{\prime }}\times [0,T] ) .  \label{A244}
\end{eqnarray}
\end{lemma}

\noindent\textbf{Proof.} Let us choose in (\ref{eq7}) the test function $%
\phi (\vec x , t):=\psi (\vec x)\, \varphi (t) $, such that $\psi (\vec x)
\in H^{1}_0(\Omega )$ and $\varphi (t) \in W^{1}_{1} ([0,T]):\;\;$ $\varphi
(0)=\varphi (T)=0$, which yields
\begin{equation*}
\Bigl|\int_0^T \Bigl( \int\limits_{\Omega }\omega \,\psi \; d\vec{x} \Bigl)%
\varphi _{t} dt\Bigl|\leqslant C\, \Vert \psi \Vert _{ H_0^{1}(\Omega ) } ||
\varphi ||_{L_1 (0,T)},
\end{equation*}
which is equivalent to (\ref{A144}).

Now we show the H\"{o}lder continuity of ${h }$ with respect to the time $%
t\in [0,T]$. By the representation (\ref{Apresent}) this can be done for $%
h_1:=(K_{1} \ast \omega)$ and $h_2:=(K_{2} \ast a)$, separately. Let
\begin{eqnarray*}
{\mathbf{\rho }} (s)\in C_0^\infty (\mathbb{R}) \;\;\;\,\;\;\mbox{ with }
\;\,\;\;\;\,\;\;\;\, {\mathbf{\rho }} (s):=\left\{
\begin{array}{ll}
1, \mbox{ if } & |s| \leqslant 1; \\
0, \mbox{ if } & |s|> 2 .%
\end{array}
\right.  \label{ddeq6sec333}
\end{eqnarray*}
We introduce the functions ${\mathbf{\rho }_{\sigma}} (\vec x):={\mathbf{%
\rho }} (|\vec x|/ \sigma )$ and ${\mathbf{\rho }}_\sigma^\Gamma (\vec x):={%
\mathbf{\rho }} (d(\vec x) / \sigma )$ for $\sigma >0$. Let $\Omega^{\prime
} $ be a sub domain of $\Omega $, such that $\delta:=dist(\overline{%
\Omega^{\prime }}, \Gamma )>0 $. Then, for any fixed $\vec x \in \overline{%
\Omega^{\prime }} \, $ and for any $\sigma <\frac \delta 4 $, the function $%
\bigtriangledown h_1$ can be written in the form
\begin{eqnarray}
\nabla h_1 (\vec{x}, t)&=&\int_{\Omega} \nabla_{\vec x}K_{1}(\vec{x},\vec{y}%
) \,\, \omega (\vec y, t)\; d\vec y =  \notag \\
&=&\int_{\Omega}\nabla_{\vec x}K_{1}(\vec{x},\vec y) \,\, {\mathbf{\rho }}%
_\sigma^\Gamma (\vec{y}) \,\, \omega (\vec y, t) \; d\vec y +  \notag \\
&+&\int_{\Omega}\nabla_{\vec x}K_{1}(\vec{x},\vec y) \,\,{\mathbf{\rho }}%
_\sigma (\vec{x}-\vec{y}) \,\, \omega (\vec y, t) \; d\vec y\,+  \notag \\
&+&\int_{\Omega}\nabla_{\vec x}K_{1}(\vec{x},\vec y) \,\, (1 - {\mathbf{\rho
}}_\sigma^\Gamma (\vec{y}) -{\mathbf{\rho }}_\sigma (\vec{x}-\vec{y})) \,\,
\omega (\vec y, t) \; d\vec y.\;\;\;\;\;\;\;\;\;\;\;  \label{eq1sec453}
\end{eqnarray}
Hence for any $t_1, t_2\in [0,T]$ we have that
\begin{equation*}
\nabla h_1 (\vec{x}, t_2)- \nabla h_1(\vec{x}, t_1)=I_1+ I_2 + I_3.
\end{equation*}
Using $\omega\in L_\infty (\Omega_T)$, the terms $I_i, \,\, i=1,2,3$ are
estimated by
\begin{eqnarray}
|I_{1}| &\leqslant& ||{\mathbf{\rho }}_\sigma^\Gamma \, \nabla_{\vec x}
K_{1} (\vec x , \cdot) ||_{L_1(\Omega )} \; ||\omega (\cdot, t_2) -\omega
(\cdot, t_1) ||_{L_\infty (\Omega )}\leqslant C\sigma ,  \notag \\
|I_{2}| &\leqslant& ||{\mathbf{\rho }_{\sigma}}\, \nabla_{\vec x} K_{1}
(\vec x , \cdot) ||_{L_1(\Omega )} \; ||\omega (\cdot, t_2) -\omega (\cdot,
t_1) ||_{L_\infty (\Omega )}\leqslant C\sigma  \label{eq1sec453}
\end{eqnarray}
and
\begin{eqnarray}
|I_3| \leqslant ||(1- {\mathbf{\rho }}_\sigma^\Gamma-{\mathbf{\rho }_{\sigma}%
}) \, \nabla_{\vec x} K_{1} (\vec x , \cdot)||_{H^1_0(\Omega )} ||\omega
(\cdot, t_2) -\omega (\cdot, t_2) ||_{H^{-1}(\Omega )} .  \label{ddd}
\end{eqnarray}
Since the function $(1-{\mathbf{\rho }}_\sigma^\Gamma-{\mathbf{\rho }%
_{\sigma}})\, \nabla_{\vec x} K_{1} (\vec x , \vec y), \,\,\, \vec y \in
\Omega $ has a singularity just at the point $\vec x \in \overline{%
\Omega^{\prime }} $, we obtain
\begin{eqnarray}
||(1-{\mathbf{\rho }}_\sigma^\Gamma-{\mathbf{\rho }_{\sigma}})\,
\nabla_{\vec x} K_{1} (\vec x , \cdot) ||_{H^1_0(\Omega )} \leqslant C \,
|\,ln(\sigma)|  \label{eq1032sec453}
\end{eqnarray}
and by (\ref{A144})
\begin{eqnarray}
|| \, \omega (\cdot, t_2) &-&\omega (\cdot, t_1) ||_{H^{-1}(\Omega )}
\leqslant ||\int_{t_1}^{t_2} \partial_t \omega (\cdot , t) \; dt
||_{H^{-1}(\Omega )}  \notag \\
&\leqslant& |t_2 -t_1 | \; \max_{t\in[0,T]} || \, \partial_t \omega \;
||_{L_\infty (0,T; \; H^{-1}(\Omega ) )} \leqslant C |t_2 -t_1 |.
\label{eq1032sec453}
\end{eqnarray}
Therefore choosing $\sigma:= |t_2 -t_1 |, $ from (\ref{eq1sec453})-(\ref%
{eq1032sec453}), we derive
\begin{eqnarray}
|\nabla h_1(\vec{x}, t_2)- \nabla h_1(\vec{x}, t_1)| \leqslant C |t_2 -t_1
|\,| \,ln |t_2 -t_1|\,|  \label{eq4sec453}
\end{eqnarray}
for all $t_1, t_2\in [0,T]$. Using this approach we can deduce the same
estimate for $h_1$.

In view (\ref{ell-6}) of Lemma \ref{teo2sec2}, we have
\begin{equation}  \label{1ell-112}
||h_2 (\cdot,t_2)- h_2 (\cdot,t_1) ||_{C^{1}(\overline{\Omega^{\prime } }%
)}\leqslant C||a(\cdot,t_2 )-a(\cdot,t_1 ) ||_{L_1(\Gamma )} \leqslant C
|t_2 -t_1 |^{\theta }
\end{equation}
for any $t_1, t_2\in [0,T]$. Therefore, accounting (\ref{eq4sec453}) and (%
\ref{1ell-112}), we derive (\ref{A244}). This completes the proof of the
Lemma. $\;\;\;\;\;\;\;\;\blacksquare$

\section{Appendix}

\label{sec6}

Let us suppose that initially the boundary and initial data $b,\,\, \omega
_{0}$ of the problem (\ref{eq1})-(\ref{eq6}) are known. In this section we
give sufficient conditions on these data under which there exists at least
one extension $\breve\omega$ on all domain $\overline\Omega_{T} $,
satisfying the conditions (\ref{eq01sec1})-(\ref{eq03sec1}).

Let us assume that the data $b,\,\, \omega_{0}$ admit the following
regularity
\begin{eqnarray}
\omega_{0}\in W_{p}^{(2-\frac{1}{p})} (\Omega)\cap L_{\infty}(\Omega),
\;\;\;\;\; b\in W_{p}^{(2-\frac{1}{p})}(\Gamma_{T}^{-}) \cap
L_{\infty}(\Gamma^{-}_{T})  \label{eq01sec6}
\end{eqnarray}
with some $p>n$, such that
\begin{equation}
\omega_{0}(\vec x)\geqslant0 \mbox{ a.e. in } \Omega,\quad b(\vec
x,t)\geqslant0 \mbox{ a.e. on } \Gamma^{-}_{T} .  \label{eq02sec6}
\end{equation}

Since $\Gamma\in C^{2+\gamma}$ for $\gamma>0$, there exists a function $%
\breve b=\breve b (\vec x ,t)$, defined on the boundary $\Gamma_{T}$, such
that
\begin{eqnarray}
\breve b\,\,\,\in & W_{p}^{(2-\frac{1}{p})}(\Gamma_{T}) \cap L_{\infty
}(\Gamma_{T}) \;\;\mbox{ and }\;\;\breve b \big|_{\Gamma_{T}^{-}} =b ;\;\;\;
\,\, \breve b \geqslant0 \,\, \mbox{ on } \,\, \Gamma_{T}  \label{eq0101sec6}
\end{eqnarray}
and
\begin{eqnarray}
\left\{
\begin{array}{lll}
||\breve b||_{L_{\infty}(\Gamma_{T})} & \leqslant & ||b||_{L_{\infty}(%
\Gamma^{-}_{T})} , \\
&  &  \\
||\breve b||_{W_{p}^{(2-\frac{1}{p})} (\Gamma_{T} )} & \leqslant &
||b||_{W_{p}^{(2-\frac{1}{p})} (\Gamma_{T}^{-} )}.%
\end{array}
\right.  \label{eq14sec6}
\end{eqnarray}
Let $\breve\omega$ be the solution of the system
\begin{eqnarray*}
\left\{
\begin{array}{ll}
\breve\omega_{t}-\Delta\breve\omega=0, \;\;\;\; (\vec x,t)\in\Omega_{T}; &
\\
&  \\
\breve\omega\big|_{\Gamma_{T}}=\breve b;\;\;\;\;\;\;\;\;\;\;\; \breve\omega%
\big|_{t=0}=\omega_{0}(\vec x), \;\;\;\vec x\in\Omega. &
\end{array}
\right.  \label{eq4sec3}
\end{eqnarray*}
According to (\ref{eq01sec6})-(\ref{eq14sec6}), by \cite{LadySolonUral68},
Theorem 9.1, p. 341 and Theorem 7.1, p.181, there exists an unique solution $%
\breve\omega\in W^{2,1}_{p} (\Omega_{T})\cap L_{\infty}(\Omega_{T})$, such
that
\begin{eqnarray}
\breve\omega &\geqslant & 0 \,\,\,\,\,\, \mbox{a.e.  on } \,\, \Omega _{T},
\notag \\
||\breve\omega||_{L_{\infty}(\Omega_{T})} & \leqslant& C(||\breve
b||_{L_{\infty}(\Gamma_{T})}+|| \omega_{0} ||_{L_{\infty}(\Omega)})  \notag
\\
& \leqslant& C(||b||_{L_{\infty}(\Gamma^{-}_{T})}+||\omega_{0}||_{L_{\infty
}(\Omega)}) \,\,\,  \label{eq55sec3}
\end{eqnarray}
and
\begin{eqnarray}
||\breve\omega||_{W^{2,1}_{p} (\Omega_{T})} & \leqslant& C( ||\breve
b||_{W_{p}^{(2-\frac{1}{p})} (\Gamma_{T} )} +||\omega_{0} ||_{W_{p}^{(2-%
\frac{2}{p})} (\Omega)} )  \notag \\
& \leqslant& C( ||b||_{W_{p}^{(2-\frac{1}{p})} (\Gamma^{-}_{T}
)}+||\omega_{0} ||_{W_{p}^{(2-\frac{2}{p})} (\Omega)} ) , \,\,\,
\label{eq5sec3}
\end{eqnarray}
where the constants $C$ depend only on $p$ and $\Omega$. Last inequality
implies
\begin{eqnarray}
\int_{0}^{T} ||\bigtriangledown\breve\omega(\cdot, t) ||^{p}_{L_{\infty
}(\Omega)} \;dt \leqslant C ||\breve\omega||^{p}_{W^{2,1}_{p} (\Omega_{T})}
. \,\,\,  \label{eq555sec3}
\end{eqnarray}
Combining the estimates (\ref{eq55sec3})-(\ref{eq555sec3}), we see that the
constructed function $\breve\omega$ is an extension of $b,\,\, \omega_{0}$
on all domain $\overline\Omega_{T} $, which satisfies the conditions (\ref%
{eq01sec1})-(\ref{eq03sec1}).

\end{document}